\documentstyle[twoside,amstex]{article}
\input amssym.def
\input amssym.tex
\def\bc{\begin{center}}
\def\ec{\end{center}}
\def\no{\noindent}

\begin{document}
\thispagestyle{empty} \vspace*{3 true cm} \pagestyle{myheadings}
\markboth {\hfill {\sl Huanyin Chen, H. Kose and Y.
Kurtulmaz}\hfill} {\hfill{\sl STRONGLY CLEAN MATRICES OVER \\
COMMUTATIVE RINGS}\hfill} \vspace*{-1.5 true cm} \bc{\large\bf
STRONGLY CLEAN MATRICES OVER\\ \vskip2mm COMMUTATIVE RINGS}\ec

\vskip6mm
\bc{{\bf Huanyin Chen}\\[1mm]
Department of Mathematics, Hangzhou Normal University\\
Hangzhou 310036, China, huanyinchen@@aliyun.com}\ec
\bc{{\bf H. Kose}\\[1mm]
Department of Mathematics, Ahi Evran University\\ Kirsehir,
Turkey, handankose@@gmail.com}\ec
\bc{{\bf Y. Kurtulmaz}\\[1mm]
Department of Mathematics, Bilkent University\\ Ankara, Turkey,
yosum@@fen.bilkent.edu.tr}\ec

\begin{abstract} A commutative ring $R$ is projective free provided that every finitely generated $R$-module is free.
An element in a ring is strongly clean provided that it is the sum
of an idempotent and a unit that commutates. Let $R$ be a
projective-free ring, and let $h\in R[t]$ be a monic polynomial of
degree $n$. We prove, in this article, that every $\varphi\in
M_n(R)$ with characteristic polynomial $h$ is strongly clean, if
and only if the companion matrix $C_h$ of $h$ is strongly clean,
if and only if there exists a factorization $h=h_0h_1$ such that
$h_0\in {\Bbb S}_0, h_1\in {\Bbb S}_1$ and $(h_0,h_1)=1$. Matrices
over power series over projective rings are also discussed. These
extend the known results [1, Theorem 12] and [5, Theorem 25].\\[2mm]
\end{abstract}

\section{Introduction}

\vskip4mm Let $R$ be a ring with an identity. An element $a\in R$
is strongly clean provided that it is the sum of an idempotent and
a unit that commutates. It is attractive to express a matrix over
a commutative ring as the sum of an idempotent matrix and an
invertible matrix that commutates. A ring $R$ is local if it has
only one maximal right ideal. As is well known, a ring $R$ is
local if and only if for any $a\in R$, either $a$ or $1-a$ is
invertible. In [9, Example 1], Wang and Chen constructed $2\times
2$ matrices over a commutative local ring which are not strongly
clean, e.g., $\left( \begin{array}{cc} 8&6\\
3&7 \end{array} \right)\in M_{(2)}(R)$. In fact, it is hard to
determine when a matrix is strongly clean. In [4, Theorem 8], Chen
et al. discussed when every $2\times 2$ matrix over a commutative
local ring $R$, i.e., $M_n(R)$, is strongly clean. In [8, Theorem
2.6], Li investigated when a single $2\times 2$ matrix over a
commutative local ring is strongly clean. In [13, Theorem 7], Yang
and Zhou characterized a $2\times 2$ matrix ring over a local ring
(not necessarily commutative) in which every matrix is strongly
clean. Strongly clean generalized $2\times 2$ matrices over a
local ring were also studied by Tang and Zhou (cf. [10, Theorem
15]. In [1, Theorem 12], Borooah et al. characterized when an
$n\times n$ matrix over a commutative local ring $R$ is strongly
clean, in terms of factorization in the polynomial ring $R[t]$. A
commutative ring $R$ is projective free provided that every
 finitely generated $R$-module is free. The class of projective
free rings is very large. For instances, commutative local rings,
principal ideal domains, B\'{e}zout domains, etc. So as to study
strong cleanness of matrices over a projective free ring, Fan and
Yang introduced the condition UCP, and characterize strong
cleanness of matrices over commutative projective free rings
having UCP (cf. [5, Theorem 25]) (We note that projective free
condition was omitted there, but it is necessary from the
procedure of its proof). Further, they asked that whether every
commutative projective free ring has ULP. We give, in this
article, an affirmative answer to this question, by a different
route. That is, we prove that every $\varphi\in M_n(R)$ with
characteristic polynomial $h$ is strongly clean, if and only if
the companion matrix $C_h$ of $h$ is strongly clean, if and only
if there exists a factorization $h=h_0h_1$ such that $h_0\in {\Bbb
S}_0, h_1\in {\Bbb S}_1$ and $(h_0,h_1)=1$. Matrices over power
series over projective rings are also discussed. These extend the
known results [1, Theorem 12] and [5, Theorem 25].

Let $h(t)\in R[t]$. We say that $h(t)$ is a monic polynomial of
degree $n$ if $h(t)=t^n+a_{n-1}t^{n-1}+\cdots +a_1t+a_0$ where
$a_{n-1},\cdots ,a_1,a_0\in R$. Every square matrix $\varphi\in
M_n(R)$ over a commutative ring $R$ is associated with a
characteristic polynomial $\chi (\varphi)$. Let $f,g\in R[t]$. The
notation $(f,g)=1$ means that there exist some $h,k\in R[t]$ such
that $fh+gk=1$. That is, the ideal generated by $f,g$ is $R[t]$.
We write $U(R)$ for the set of all invertible elements in $R$ and
$M_n(R)$ for the rings of all $n\times n$ matrices over a ring
$R$. $R[t]$ and $R[[t]]$ always stands for the rings of
polynomials and power series over a ring $R$, respectively.

\section{Strong Cleanness in $M_n(R)$}

\vskip4mm Let $h(t)=t^n+a_{n-1}t^{n-1}+\cdots +a_1t+a_0$. We call
the matrix
$$\varphi =\left(
\begin{array}{ccccc}
0&0&\cdots &0&-a_{0}\\
1&0&\cdots &0&-a_{1}\\
\vdots&\vdots&\ddots &\vdots&\vdots\\
0&0&\cdots &1&-a_{n-1}
\end{array}
\right)\in M_n(R)$$ is the companion matrix of $h$, and denote it
by $C_h$. The minimal polynomial $p(t)\in F[t]$ of a matrix
$\varphi$ over the filed $F$ is the monic polynomial in $\varphi$
of smallest degree such that $p(\varphi)=0$. We note it by
$m(\varphi)$. The characteristic polynomial as well as the minimal
polynomial of the companion matrix $C_h$ over a field $F$ are
equal to $h$. In this section we completely determine a strongly
clean matrix over a projective free ring, in terms of
factorization of its characteristic polynomial.

Let $R$ be a commutative ring. Given polynomials
$f(t)=t^m+a_1t^{m-1}+\cdots +a_m, g(t)=b_0t^n+b_1t^{n-1}+\cdots
+b_n (b_0\neq 0)\in R[t]$, the resultant of $f$ and $g$ is defined
by the determinant of the $(m+n)\times (m+n)$ matrix
$$res(f,g)=\left|
\begin{array}{ccccccccc}
1&a_1&\cdots &a_m&&&&\\
&1&a_1&\cdots &a_m&&&\\
&&1&a_1&\cdots &a_m&&\\
&&&\ddots&\ddots&\ddots &\ddots&&\\
&&&&1&a_1&\cdots&&a_m\\
b_0&b_1&\cdots &b_n&&&&\\
&b_0&b_1&\cdots &b_n&&&\\
&&b_0&b_1&\cdots &b_n&&\\
&&&\ddots&\ddots&\ddots &\ddots&&\\
&&&&b_0&b_1&\cdots&&b_n\\
\end{array}
\right|$$ where blank spaces consist of zeros. The following is
called the Weyl Principal. Let $f,g$ be polynomials in ${\Bbb
Z}[x_1,x_2,\cdots ,x_n]$ with $g\neq 0$ (i.e., $g$ is not the zero
polynomial). If for all $(\alpha_1,\cdots ,\alpha_n)\in {\Bbb
Q}^n$, $f(\alpha_1,\cdots ,\alpha_n)=0$ whenever
$g(\alpha_1,\cdots ,\alpha_n)\neq 0$, then $f=0$ in any
commutative ring $R$. We begin with the following results which
are analogous to those over fields.

\vskip4mm \hspace{-1.8em} {\bf Lemma 2.1.}\ \ {\it Let $R$ be a
commutative ring, and let $f\in R[t]$ be monic and $g,h\in R[t]$.
Then the following are equivalent:}\vspace{-.5mm}
\begin{enumerate}
\item [(1)]{\it $res(f,g)=res(f,g+fh)$.}
\vspace{-.5mm}
\item [(2)]{\it $res(f,gh)=res(f,g)res(f,h)$.}
\end{enumerate}\vspace{-.5mm} {\it Proof.}\ \ $(1)$ Write $h=c_0t^s+\cdots +c_s\in R[t]$. It will suffice to
show that $res(f,g)=res(f,g+c_{s-i}t^if)$. Since any determinant
in which every entry in a row is a sum of two elements is the sum
of two corresponding determinants, the result follows.

$(2)$ Write $f=t^m+a_1t^{m-1}+\cdots +a_m,
g=b_0t^n+b_1t^{n-1}+\cdots +b_n, h=c_0t^s+c_1t^{s-1}+\cdots +c_s$.
Then
$$\begin{array}{lll}
\alpha (a_1,\cdots ,a_m;b_0,\cdots ,b_n;c_0,\cdots ,c_s):&=&res(f,gh)-res(f,g)res(f,h)\\
&\in &{\Bbb Z}[a_1,\cdots ,a_m;b_0,\cdots ,b_n;c_0,\cdots ,c_s].
\end{array}$$
Consider $\alpha (x_1,\cdots ,x_m;y_0,\cdots ,y_n;z_0,\cdots ,z_s)
\in {\Bbb Z}[x_1,\cdots ,x_m;y_0,\cdots ,y_n;$ $z_0,\cdots ,z_s]$.
Clearly, the result holds if $R={\Bbb Q}$.

For any $u_1,$ $\cdots ,u_m;$ $v_0,$ $\cdots ,v_n;$ $w_0,\cdots
,w_s\in {\Bbb Q}$, we see that $$\alpha (u_1,\cdots
,u_m;v_0,\cdots ,v_n;w_0,\cdots ,w_s)=0.$$ By the Weyl Principal
that $$\alpha (x_1,\cdots ,x_m;y_0,\cdots ,y_n;z_0,\cdots
,z_s)=0$$ in ${\Bbb Z}[x_1,\cdots ,x_m;y_0,\cdots ,y_n;z_0,\cdots
,z_s]$. Therefore $\alpha (a_1,\cdots ,a_m;b_0,\cdots ,b_n;$ $
c_0,\cdots ,c_s)$ $=0$, and so
$res(f,gh)=res(f,g)res(f,h)$.\hfill$\Box$

\vskip4mm \hspace{-1.8em} {\bf Lemma 2.2.}\ \ {\it Let $R$ be a
commutative ring, and let $f,g\in R[t]$ be monic. Then the
following are equivalent:}\vspace{-.5mm}
\begin{enumerate}
\item [(1)]{\it $(f,g)=1$.}
\vspace{-.5mm}
\item [(2)]{\it $res(f,g)\in U(R)$.}
\end{enumerate}\vspace{-.5mm} {\it Proof.}\ \ $(1)\Rightarrow (2)$ As $(f,g)=1$, we can
find some $u,v\in R[t]$ such that $uf+vg=1$. By virtue of Lemma
2.1, one easily checks that
$res(f,vg)=res(f,v)res(f,g)=res(f,vg+uf)=res(f,1)=1.$ Accordingly,
$res(f,g)\in U(R)$.

$(2)\Rightarrow (1)$ Let $m=deg(f)$ and $n=deg(g)$. Observing that
$$res(f,g)
\left|
\begin{array}{cc}
I_{m+n-1}&
\begin{array}{c}
t^{m+n}\\
\vdots\\
t
\end{array}\\
\begin{array}{ccc}
0&\cdots&0
\end{array}&1
\end{array}
\right|=\left|
\begin{array}{cc}
*&
\begin{array}{c}
t^{n}f\\
\vdots\\
f\\
t^{m}g\\
\vdots\\
g
\end{array}
\end{array}
\right|,$$ therefore we can find some $u,v\in R[t]$ such that
$res(f,g)=uf+vg$. Hence
$\big(res(f,g)\big)^{-1}uf+\big(res(f,g)\big)^{-1}vg=1$, as
asserted.\hfill$\Box$

\vskip4mm \hspace{-1.8em} {\bf Lemma 2.3.}\ \ {\it Let $R$ be a
commutative ring, and let $f,g\in R[t]$ be monic. Then the
following are equivalent:}\vspace{-.5mm}
\begin{enumerate}
\item [(1)]{\it $(f,g)=1$.}
\vspace{-.5mm}
\item [(2)]{\it For all maximal ideal $M$ of $R$, $\big(\overline{f},\overline{g}\big)=\overline{1}$.}
\end{enumerate}\vspace{-.5mm} {\it Proof.}\ \ $(1)\Rightarrow (2)$
is obvious.

$(2)\Rightarrow (1)$ If $(f,g)\neq 1$, it follows from Lemma 2.2
that $res(f,g)\not\in U(R)$. Thus, $res(f,g)R\neq R$, and so we
can find a maximal ideal $M$ of $R$ such that $res(f,g)R\subseteq
M\subsetneq R$. Hence, $\overline{res(f,g)}=\overline{0}$ in
$R/M$. This implies that
$ref\big(\overline{f},\overline{g}\big)=\overline{0}$. As $R/M$ is
a filed, by using Lemma 2.2 again,
$\big(\overline{f},\overline{g}\big)\neq \overline{1}$, a
contradiction, and thus yielding the result.\hfill$\Box$

\vskip4mm For $r\in R$, define
$${\Bbb S}_r=\{ f\in R[t]~|~f~\mbox{monic, and}~f(r)\in U(R)~\}.$$

\vskip4mm \hspace{-1.8em} {\bf Theorem 2.4.}\ \ {\it Let $R$ be a
projective-free ring, and let $h\in R[t]$ be a monic polynomial of
degree $n$. Then the following are equivalent:} \vspace{-.5mm}
\begin{enumerate}
\item [(1)]{\it Every $\varphi\in M_n(R)$ with $\chi (\varphi)=h$ is strongly clean.}
\vspace{-.5mm}
\item [(2)]{\it The companion matrix $C_h$ of $h$ is strongly clean.}\vspace{-.5mm}
\item [(3)]{\it There exists a factorization $h=h_0h_1$ such that $h_0\in {\Bbb
S}_0, h_1\in {\Bbb S}_1$ and $(h_0,h_1)=1$.}
\end{enumerate}
\vspace{-.5mm} {\it Proof.}\ \ $(1)\Rightarrow (2)$ Write
$h=t^n+a_{n-1}t^{n-1}+\cdots +a_1t+a_0\in R[t]$. Then
$$C_h=\left(
\begin{array}{ccccc}
0&0&\cdots &0&-a_{0}\\
1&0&\cdots &0&-a_{1}\\
\vdots&\vdots&\ddots &\vdots&\vdots\\
0&0&\cdots &1&-a_{n-1}
\end{array}
\right)\in M_n(R).$$ Clearly, $\chi (C_h)=h$. Moreover,
$\varphi\in M_n(R)$ is strongly clean by hypothesis.

$(2)\Rightarrow (3)$ In view of [1, Lemma 1], there exists a
decomposition $nR=A\oplus B$ such that $A$ and $B$ are
$C_h$-invariant, $C_h~|_{A}\in aut(A)$ and $(I_n-C_h)~|_{B}\in
aut(B)$. Since $R$ is projective-free, there exist $p,q\in {\Bbb
N}$ such that $A\cong pR$ and $B\cong qR$. As $R$ is commutative,
we see that $p+q=n$. Regarding $end_R(A)$ as $M_p(R)$, we see that
$C_h~|_{A}\in GL_p(R)$. In light of [1, Lemma 4], $\chi
(C_h~|_{A})\in {\Bbb S}_0$. Similarly, $(I_n-C_h)~|_{B}\in
GL_q(R)$. This implies that $det\big((I_n-C_h)~|_{B}\big)\in
U(R)$, and so $det\big(t I_q-C_h~|_{B}\big)\in {\Bbb S}_1$. Hence,
we get $\chi (C_h~|_{B})\in {\Bbb S}_1$. Clearly, $\chi (C_h)=\chi
(C_h~|_{A})\chi (C_h~|_{B})$. Choose $h_0=\chi (C_h~|_{A})$ and
$h_1=\chi (C_h~|_{B})$. Then $h=h_0h_1$. Let $M$ be a maximal
ideal of the ring $R$, and let $F=R/M$. Assume that
$\big(\overline{h_0},\overline{h_1}\big)\neq \overline{1}$ in
$F[t]$. Clearly, $F[t]$ is a principal ideal domain. Write
$\big(\overline{h_0},\overline{h_1}\big)=\big(p(t)\big)$, where
$p(t)\in F[t]$ is monic and $degp(t)\geq 1$. As $F$ is a field,
$m\big(\overline{C_h}\big)=m\big(C_{\overline{h}}\big)=\chi\big(C_{\overline{h}}\big)=\overline{h}$.
As $A$ and $B$ are $C_h$-invariant respectively, we see that $C_h$
is similar to $\left(
\begin{array}{cc}
C_h~|_{A}&\\
&C_h~|_{B} \end{array} \right)$. Thus, we get $\overline{C_h}$ is
similar to $\left(
\begin{array}{cc}
\overline{C_h~|_{A}}&\\
&\overline{C_h~|_{B}}\end{array} \right)$. Accordingly, we get
$m\big(C_{\overline{h}}\big)=\big[m(\overline{C_h~|_{A}}),m(\overline{C_h~|_{B}})\big]$.
But we also have
$\chi\big(\overline{C_h~|_{A}}\big)=\overline{\big(\chi
C_h~|_{A}\big)}=\overline{h_0}$. Likewise,
$\chi\big(\overline{C_h~|_{B}}\big)=\overline{h_1}$. Hence,
$\overline{h}|~\big[\overline{h_0},\overline{h_1}\big]$. As
$p(t)~|~\overline{h_0},\overline{h_1}$, we see that
$deg(h)=deg(\overline{h})\leq
deg\big[\overline{h_0},\overline{h_1}\big] <
deg\big(\overline{h_0}\overline{h_1}\big)=deg\big(h_0h_1\big)=n$,
a contradiction. Therefore
$\big(\overline{h_0},\overline{h_1}\big)=\overline{1}$ for all
maximal ideals $M$ of $R$. By Lemma 2.4, we have that
$(h_0,h_1)=1$, as desired.

$(3)\Rightarrow (1)$ is clear from [5, Corollary 4].\hfill$\Box$

\vskip4mm \hspace{-1.8em} {\bf Example 2.5.}\ \ Let $A= \left(
\begin{array}{cc}
0&3\\
1&2
\end{array}
\right)\in M_2({\Bbb Z})$. Then $\chi(A)=t^2-2t-3$. It is easy to
verify that there are no any $h_0\in {\Bbb S}_0$ and $h_1\in {\Bbb
S}_1$ such that $\chi(A)=h_0h_1$. Accordingly, $A\in M_2({\Bbb
Z})$ is not strongly clean by Theorem 2.4.

\vskip4mm A matrix $\varphi\in M_n(R)$ is cyclic if there exists a
column $\alpha$ such that $(\alpha,\varphi\alpha,$ $\cdots,
\varphi^{n-1}\alpha)\in GL_n(R)$. Let $\varphi=(\varphi_{ij})\in
M_n(R)$ be a matrix which coincides with a companion matrix below
the main diagonal, i.e., $\varphi_{ij}=1, i=j+1$ and
$\varphi_{ij}=0$ when $i\geq j+2$. We note that $\varphi$ is
cyclic.

\vskip4mm \hspace{-1.8em} {\bf Corollary 2.6.}\ \ {\it Let $R$ be
a projective-free ring, and let $\varphi\in M_n(R)$ be a cyclic
matrix. Then the following are equivalent:} \vspace{-.5mm}
\begin{enumerate}
\item [(1)]{\it $\varphi\in M_n(R)$ is strongly clean.}
\vspace{-.5mm}
\item [(2)]{\it There exists a factorization $\chi(\varphi )=h_0h_1$ such that $h_0\in {\Bbb
S}_0, h_1\in {\Bbb S}_1$ and $(h_0,h_1)=1$.}
\end{enumerate}
\vspace{-.5mm} {\it Proof.}\ \ $(1)\Rightarrow (2)$ Let $F=nR$ be
free as a right $R$-module, and let $\{ e_1,\cdots ,e_n\}$ be the
standard basis of $F$. Define $\sigma :F\to F$ given by $\sigma
(e_1,\cdots ,e_n)=(e_1,\cdots ,e_n)\varphi$. Since $\varphi$ is
cyclic, there exists a column $\alpha$ such that
$(\alpha,\varphi\alpha,\cdots,$ $\varphi^{n-1}\alpha)\in GL_n(R)$.
Thus, $\{ \alpha,\varphi\alpha,\cdots, \varphi^{n-1}\alpha \}$ is
a basis of $F$. So there are some $c_1,\cdots ,c_n\in R$ such that
$\varphi^n\alpha=\alpha c_1+\varphi\alpha c_2+\cdots
+\varphi^{n-1}\alpha c_n$. Clearly, $(\alpha,\varphi\alpha,\cdots,
\varphi^{n-1}\alpha)=(e_1,\cdots
,e_n)(\alpha,\varphi\alpha,\cdots, \varphi^{n-1}\alpha)$. Thus,
$$\begin{array}{ccl}
\sigma(\alpha,A\alpha,\cdots,
\varphi^{n-1}\alpha)&=&\big(\sigma(e_1,\cdots ,e_n)\big)(\alpha,\varphi\alpha,\cdots, \varphi^{n-1}\alpha)\\
&=&(e_1,\cdots ,e_n)\varphi(\alpha,\varphi\alpha,\cdots, \varphi^{n-1}\alpha)\\
&=&(\varphi\alpha,\varphi^2\alpha,\cdots, \varphi^{n-1}\alpha,\varphi^n\alpha)\\
&=&(\alpha,\varphi\alpha,\cdots, \varphi^{n-1}\alpha)\left(
\begin{array}{ccccc}
0&0&\cdots &0&c_1\\
1&0&\cdots &0&c_2\\
0&1&\cdots &0&c_3\\
\vdots&\vdots &\ddots &\vdots&\vdots\\
0&0&\cdots &1&c_n
\end{array}
\right).
\end{array}$$
Let $\gamma=(\alpha,\varphi\alpha,\cdots,$
$\varphi^{n-1}\alpha)\in GL_n(R)$. Then
$$\gamma^{-1}\varphi\gamma=\left(
\begin{array}{ccccc}
0&0&\cdots &0&c_1\\
1&0&\cdots &0&c_2\\
0&1&\cdots &0&c_3\\
\vdots&\vdots &\ddots &\vdots&\vdots\\
0&0&\cdots &1&c_n
\end{array}
\right).$$ As $\varphi\in M_n(R)$ is strongly clean, it follows
that the preceding companion matrix of $\chi (\varphi)$ is
strongly clean. This substantiates our claim by Theorem 2.4.

$(2)\Rightarrow (1)$ is trivial by Theorem 2.4.\hfill$\Box$

\section{Matrices over Power Series}

\vskip4mm The purpose of this section is to investigate strong
cleanness of $n\times n$ matrices over power series over a
commutative projective free rings. Let
$A(x)=\big(a_{ij}(x)\big)\in M_n\big(R[[x]]\big)$, where each
$a_{ij}(x)\in R[[x]]$. We use $A(0)$ to denote the matrix
$\big(a_{ij}(0)\big)\in M_n(R)$. Then we have

\vskip4mm \hspace{-1.8em} {\bf Theorem 3.1.}\ \ {\it Let $R$ be a
projective-free ring, and let $A(x)\in M_n\big(R[[x]]\big) (n\in
{\Bbb N})$. Then the following are equivalent:} \vspace{-.5mm}
\begin{enumerate}
\item [(1)]{\it $A(0)\in M_n(R)$ is strongly clean.} \vspace{-.5mm} \item
[(2)]{\it $A(x)\in M_n\big(R[[x]]\big)$ is strongly clean.}
\end{enumerate}
\vspace{-.5mm}  {\it Proof.}\ \ $(1)\Rightarrow (2)$ Obviously,
$R[[x]]$ is projective-free. Let $H(x,t)=\chi \big(A(x)\big)\in
R[[x]][t]$. Then $H(0,t)=\chi \big(A(0)\big)\in R[t]$. By using
Theorem 2.4, $H(0,t)=h_0h_1$, where
$h_0=t^m+\alpha_1t^{m-1}+\cdots +\alpha_m \in {\Bbb S}_0,
h_1=t^s+\beta_1t^{s-1}+\cdots +\beta_s\in {\Bbb S}_1$ and
$(h_0,h_1)=1$. Next, we will find a factorization $H(x,t)=H_0H_1$
where $H_0(x,t)=t^m+\sum\limits_{i=0}^{m-1}A_i(x)t^i\in {\Bbb
S}_0$ and $H_1(x,t)= t^{s}+\sum\limits_{i=0}^{s-1}B_i(x)t^i\in
{\Bbb S}_1$. Choose $H_0(0,t)\equiv h_0$ and $H_1(0,t)\equiv h_1$.
Write
$H(x,t)=\sum\limits_{i=0}^{n}\big(\sum\limits_{j=0}^{\infty}c_{ij}x^j\big)t^i$.
Then $$\begin{array}{lll}
H(x,t)&=&\sum\limits_{j=0}^{\infty}\big(\sum\limits_{i=0}^{n}c_{ij}t^i\big)x^j\\
&=&H(0,t)+\sum\limits_{j=1}^{\infty}\big(\sum\limits_{i=0}^{n}c_{ij}t^i\big)x^j.
\end{array}$$
Write $A_i(x)=\sum\limits_{j=0}^{\infty}a_{ij}x^j$ and
$B_i(x)=\sum\limits_{j=0}^{\infty}b_{ij}x^j$. Then
$$\begin{array}{lll}
H_0&=&t^m+\sum\limits_{i=0}^{m-1}\big(\sum\limits_{j=0}^{\infty}a_{ij}x^j\big)t^i\\
&=&t^m+\sum\limits_{j=0}^{\infty}\big(\sum\limits_{i=0}^{m-1}a_{ij}t^i\big)x^j\\
&=&h_0+\sum\limits_{j=1}^{\infty}\big(\sum\limits_{i=0}^{m-1}a_{ij}t^i\big)x^j.
\end{array}$$ Likewise,
$$H_1=h_1+\sum\limits_{j=1}^{\infty}\big(\sum\limits_{i=0}^{s-1}b_{ij}t^i\big)x^j.$$
Write $H_0H_1=h_0h_1+\sum\limits_{j=1}^{\infty}z_jx^j$. Thus, we
should have
$$\begin{array}{lll}
z_1&=&h_0\big(\sum\limits_{i=0}^{s-1}b_{i1}t^i\big)+\big(\sum\limits_{i=0}^{m-1}a_{i1}t^i\big)h_1\\
&=&\sum\limits_{i=0}^{n-1}c_{i1}t^i.
\end{array}$$
This implies that
$$\big(b_{(s-1)1},\cdots ,b_{01},a_{(m-1)1},\cdots ,a_{01}\big)A=\big(c_{(n-1)1},\cdots ,c_{01}\big),$$
where $A= \left(
\begin{array}{ccccccc}
1&\alpha_1&\cdots&\alpha_m&&&\\
&1&\alpha_1&\cdots&\alpha_m&&\\
&&\ddots&\ddots&\ddots&\ddots&\\
&&&1&\alpha_1&\cdots&\alpha_m\\
1&\beta_1&\cdots&\beta_s&&&\\
&1&\beta_1&\cdots&\beta_s&&\\
&&\ddots&\ddots&\ddots&\ddots&\\
&&&1&\beta_1&\cdots&\beta_s
\end{array}
\right)$. As $(h_0,h_1)=1$, it follows from Lemma 2.2 that
$res(h_0,h_1)\in U(R)$. Thus, $det(A)\in U(R)$, and so we can find
$a_{i1},b_{j1}\in R$.
$$\begin{array}{lll}
z_2&=&h_0\big(\sum\limits_{i=0}^{s-1}b_{i2}t^i\big)+\big(\sum\limits_{i=0}^{m-1}a_{i1}t^i\big)
\big(\sum\limits_{i=0}^{s-1}b_{i1}t^i\big)+
\big(\sum\limits_{i=0}^{m-1}a_{i2}t^i\big)h_1\\
&=&\sum\limits_{i=0}^{n-1}c_{i2}t^i.
\end{array}$$
Hence $$\begin{array}{lll}
h_0\big(\sum\limits_{i=0}^{s-1}b_{i2}t^i\big)+\big(\sum\limits_{i=0}^{m-1}a_{i2}t^i\big)h_1&=&
\sum\limits_{i=0}^{n-1}c_{i2}t^i-\big(\sum\limits_{i=0}^{m-1}a_{i1}t^i\big)
\big(\sum\limits_{i=0}^{s-1}b_{i1}t^i\big)\\
&=&\sum\limits_{i=0}^{n-1}d_{i2}t^i.
\end{array}$$ Thus,
$$\big(\big(b_{(s-1)2},\cdots ,b_{02},a_{(m-1)2},\cdots ,a_{02}\big)\big)A=\big(d_{(n-1)2},\cdots ,d_{02}\big),$$
whence we can find $a_{i2},b_{j2}\in R$. By iteration of this
process, we can find $a_{ij},b_{ij}\in R, j=3,4,\cdots $.
Therefore we have $H_0$ and $H_1$ such that $H(x,t)=H_0H_1$.
Further, $H_0(x,0)=H_0(0,0)+xf(x)=h_0(0)+xf(x)\in
U\big(R[[x]]\big)$ and $H_1(x,1)=H_1(0,1)+xg(x)=h_1(1)+xg(x)\in
U\big(R[[x]]\big)$. Thus, $H_0(x,t)\in {\Bbb S}_0$ and
$H_1(x,t)\in {\Bbb S}_1$. As $(h_0,h_1)=1$, we get
$\big(H_0,H_1\big)\equiv 1 \big(mod~(xR[[x]])[t]\big)$, and so
$\big(H_0,H_1\big)+J\big(R[[x]]\big)R[[[x]][t]=R[[x]][t]$. Set
$M=R[[x]][t]/\big(H_0,H_1\big)$. Then $M$ is a finitely generated
$R[[x]]$-module, and that $J\big(R[[x]]\big)M=M$. By Nakayama's
Lemma, $M=0$, and so $\big(H_0,H_1\big)=1$. Accordingly, $A(x)\in
M_n\big(R[[x]]\big)$ is strongly clean by Theorem 2.4.

$(2)\Rightarrow (1)$ is obvious. \hfill$\Box$

\vskip4mm \hspace{-1.8em} {\bf Corollary 3.2.}\ \ {\it Let $R$ be
a projective-free ring, and let $A(x)\in M_n\big(R[x]/(x^n)\big)
(n\in {\Bbb N})$. Then the following are equivalent:}
\vspace{-.5mm}
\begin{enumerate}
\item [(1)]{\it $A(0)\in M_n(R)$ is strongly clean.} \vspace{-.5mm} \item
[(2)]{\it $A(x)\in M_n\big(R[x]/(x^n)\big)$ is strongly clean.}
\end{enumerate}
\vspace{-.5mm}  {\it Proof.}\ \ $(1)\Rightarrow (2)$ Write
$A(x)=\sum\limits_{i=0}^{\infty}a_ix^i\in
M_n\big(R[x]/(x^n)\big)$. Then $A(x)\in M_n\big(R[[x]]\big)$. In
view of Theorem 3.1, there exist
$E^2=E=\big(\sum\limits_{k=0}^{\infty}e^{ij}_kx^k\big),
U=\big(\sum\limits_{k=0}^{\infty}u^{ij}_kx^k\big)\in
GL_n\big(R[[x]]\big)$ such that $A(x)=E+U$ and $EU=UE$. As
$R[[x]]/(x^n)\cong R[x]/(x^n)$, we see that
$A(x)=\overline{E}+\overline{U}$ and
$\overline{EU}=\overline{UE}$, where
$\overline{E}^2=\overline{E}=\big(\sum\limits_{k=0}^{n-1}e^{ij}_kx^k\big)\in
M_n\big(R[x]/(x^n)\big)$ and
$\overline{U}=\big(\sum\limits_{k=0}^{n-1}u^{ij}_kx^k\big)\in
GL_n\big(R[[x]]/(x^n)\big)$, as desired.

$(2)\Rightarrow (1)$ is clear.\hfill$\Box$

\vskip4mm We now extend [6, Theorem 2.10] and [10, Theorem 2.7] as
follows.

\vskip4mm \hspace{-1.8em} {\bf Corollary 3.3.}\ \ {\it Let $R$ be
a projective-free ring, and let $n\in {\Bbb N}$. Then the
following are equivalent:} \vspace{-.5mm}
\begin{enumerate}
\item [(1)]{\it $M_n(R)$ is strongly clean.} \vspace{-.5mm} \item
[(2)]{\it $M_n\big(R[[x]]\big)$ is strongly clean.}
 \vspace{-.5mm} \item
[(3)]{\it $M_n\big(R[x]/(x^m)\big) (m\in {\Bbb N}$ is strongly
clean.}
 \vspace{-.5mm} \item
[(3)]{\it $M_n\big(R[[x_1,\cdots ,x_m]]\big) (m\in {\Bbb N}$ is
strongly clean.}
 \vspace{-.5mm} \item
[(3)]{\it $M_n\big(R[[x_1,\cdots ,x_m]]/(x_1^{n_1},\cdots
,x_m^{n_m})\big) (m\in {\Bbb N})$ is strongly clean.}
\end{enumerate}
\vspace{-.5mm}  {\it Proof.}\ \ These are obvious by induction,
Theorem 3.1 and Corollary 3.2.\hfill$\Box$

\vskip4mm \hspace{-1.8em} {\bf Example 3.4.}\ \ {\it Let $A(x)\in
M_2\big({\Bbb Z}[[x]]\big)$. Then $A(x)\in M_2\big({\Bbb
Z}[[x]]\big)$ is strongly clean if and only if $A(0)\in GL_2({\Bbb
Z})$, or $I_2-A(0)\in GL_2({\Bbb Z})$, or $A(0)$ is similar to one
of the matrices in the set $\big\{ \left(
\begin{array}{cc}
0&0\\
0&1
\end{array}
\right), \left(
\begin{array}{cc}
0&0\\
0&-1
\end{array}
\right), \left(
\begin{array}{cc}
2&0\\
0&1
\end{array}
\right), \left(
\begin{array}{cc}
2&0\\
0&-1
\end{array}
\right)\big\}.$} \vskip2mm\hspace{-1.8em} {\it Proof.}\ \ Clearly,
${\Bbb Z}$ is a principal ideal domain. It is a projective free
ring. In light of Theorem 3.1, $A(x)\in M_2\big({\Bbb
Z}[[x]]\big)$ is strongly clean if and only if so is $A(0)$.
Therefore we complete the proof, by [2, Example
16.4.9].\hfill$\Box$

\vskip4mm \hspace{-1.8em} {\bf Lemma 3.5.}\ \ {\it Let $R$ be a
commutative projective free ring, $char(R)=2$, and let $G=\{
1,g\}$ be a group. Then the following hold:} \vspace{-.5mm}
\begin{enumerate}
\item [(1)]{\it $R[x]/(x^2-1)\cong RG$.} \vspace{-.5mm}
\item [(2)]{\it $a+bg\in U(RG)$ if and only if $a+b\in U(R)$.}\vspace{-.5mm}
\item [(3)]{\it $RG$ is a projective free ring.} \vspace{-.5mm}
\end{enumerate}\vspace{-.5mm} {\it Proof.}\ \ $(1)$ is proved in [3, Lemma 2.1].

$(2)$ Obviously, $(a+bg)(a-bg)=a^2-b^2=(a+b)(a-b)$. Hence,
$(a+bg)^2=(a+b)^2$, as $char(R)=2$. If $a+bg\in U(RG)$, then
$(a+bg)(x+yg)=1$ for some $x,y\in R$. This implies that
$(a+bg)^2(x+yg)^2=1$, hence that $(a+b)^2(x+y)^2=1$. Accordingly,
$a+b\in U(R)$. The converse is analogous.\hfill$\Box$

$(3)$ Construct a ring morphism $\varphi: RG\rightarrow R$
$$\hspace{8mm}\sum\limits_{g\in G}r_gg\mapsto \sum\limits_{g\in G}r_g.$$
Then $ker(\varphi)=\{
\sum\limits_{g}r_{g}g~|~\sum\limits_{g}r_{g}=0\}$ is an ideal of
$RG$. If $a+bg\in ker(\omega)$, then $a+b=0$. Hence,
$(a+bg)^2=a^2+b^2=(a+b)^2=0$. This implies that $R\cong
RG/ker(\omega)$ with $ker(\omega)\subseteq J(RG)$. Consequently,
$RG$ is projective-free.\hfill$\Box$

\vskip4mm Let $A(x)=\big(\overline{a_{ij}(x)}\big)\in
M_n\big(R[x]/(x^2-1)\big)$ where $deg\big(a_{ij}(x)\big)\leq 1$,
and let $r\in R$. We use $A(r)$ to stand for the matrix
$\big(a_{ij}(r)\big)\in M_n(R)$.

\vskip4mm \hspace{-1.8em} {\bf Theorem 3.6.}\ \ {\it Let $R$ be a
projective-free ring with $char(R)=2$, and let $A(x)\in
M_n\big(R[[x]]/(x^2-1)\big) (n\in {\Bbb N})$. Then the following
are equivalent:} \vspace{-.5mm}
\begin{enumerate}
\item [(1)]{\it $A(1)\in M_n(R)$ is strongly clean.} \vspace{-.5mm} \item
[(2)]{\it $A(x)\in M_n\big(R[[x]]/(x^2-1)\big)$ is strongly
clean.}
\end{enumerate}
\vspace{-.5mm}  {\it Proof.}\ \ $(1)\Rightarrow (2)$ In view of
Lemma 3.5, $R[[x]]/(x^2-1)\cong RG$ is projective-free, where
$G=\{ 1,g\}$. Let $H(g,t)=\chi \big(A(g)\big)\in (RG)[t]$. Then
$H(1,t)=\chi \big(A(1)\big)\in R[t]$. In light of Theorem 2.4,
$H(1,t)=h_0h_1$, where $h_0=t^m+\alpha_{m-1}t^{m-1}+\cdots
+\alpha_0 \in {\Bbb S}_0, h_1=t^s+\beta_{s-1}t^{s-1}+\cdots
+\beta_0\in {\Bbb S}_1$ and $(h_0,h_1)=1$. We shall find a
factorization $H(g,t)=H_0H_1$ where
$H_0(g,t)=t^m+\sum\limits_{i=0}^{m-1}\big(y_i+(\alpha_i-y_i)g\big)t^i\in
{\Bbb S}_0$ and $H_1(g,t)=
t^{s}+\sum\limits_{i=0}^{s-1}\big(z_i+(\beta_i-z_i)g\big)t^i\in
{\Bbb S}_1$. Clearly, $H_0(1,t)\equiv h_0$ and $H_1(1,t)\equiv
h_1$. We will suffice to find $y_i's$ and $z_i's$. Write
$H(g,t)=\sum\limits_{i=0}^{n}\big(r_i+s_ig\big)t^i$. The equality
$H(g,t)=H_0H_1$ is equivalent to
$$\begin{array}{l}
t^n+\sum\limits_{i=0}^{n-1}r_it^i=\big(t^m+\sum\limits_{i=0}^{m-1}y_it^i\big)\big(t^s+\sum\limits_{i=0}^{s-1}z_it^i\big)+
\big(\sum\limits_{i=0}^{m-1}(\alpha_i-y_i)t^i\big)\big(\sum\limits_{i=0}^{s-1}(\beta_i-z_i)t^i\big)(*)\\
\sum\limits_{i=0}^{n-1}s_it^i=\big(t^m+\sum\limits_{i=0}^{m-1}y_it^i\big)\big(\sum\limits_{i=0}^{s-1}(\beta_i-z_i)t^i\big)+
\big(t^s+\sum\limits_{i=0}^{s-1}z_it^i\big)\big(\sum\limits_{i=0}^{m-1}(\alpha_i-y_i)t^i\big)(**).
\end{array}$$
$(**)$ holds from $H(1,t)=h_0h_1=H_0(1,t)H_1(1,t)$. $(*)$ is
equivalent to
$$\begin{array}{c}
y_0z_0+(\alpha_0-y_0)(\beta_0-z_0)=r_0,\\
y_0z_1+y_1z_0+(\alpha_0-y_0)(\beta_1-z_1)+(\alpha_1-y_1)(\beta_0-z_0)=r_1,\\
\vdots\\
y_{m-2}+y_{m-1}z_{s-1}+z_{s-2}+(\alpha_{m-1}-y_{m-1})(\beta_{s-1}-z_{s-1})=r_{n-2},\\
y_{m-1}+z_{s-1}=r_{n-1}.
\end{array}$$
As $char(R)=2$, we have
$$\begin{array}{c}
\beta_0y_0+\alpha_0z_0=r_0+\alpha_0\beta_0,\\
\beta_0y_1+\beta_1y_0+\alpha_0z_1+\alpha_1z_0=r_1+\alpha_0\beta_1+\alpha_1\beta_0,\\
\vdots\\
\beta_{s-1}y_{m-1}+y_{m-2}+\alpha_{m-1}z_{s-1}+z_{s-2}=r_{n-2}+\alpha_{m-1}\beta_{s-1},\\
y_{m-1}+z_{s-1}=r_{n-1}. \end{array}$$ This implies that
$$\big(y_{m-1},\cdots ,y_{0},z_{s-1},\cdots ,z_{0}\big)A=\big(*,\cdots ,*\big),$$
where $$A=\left(
\begin{array}{ccccccc}
1&\beta_{m-1}&\cdots&\beta_0&&&\\
&1&\beta_{m-1}&\cdots&\beta_0&&\\
&&\ddots&\ddots&\ddots&\ddots&\\
&&&1&\beta_{m-1}&\cdots&\beta_0\\
1&\alpha_{s-1}&\cdots&\alpha_0&&&\\
&1&\alpha_{s-1}&\cdots&\alpha_0&&\\
&&\ddots&\ddots&\ddots&\ddots&\\
&&&1&\alpha_{s-1}&\cdots&\alpha_0
\end{array}
\right).$$ As $(h_1,h_0)=1$, it follows from Lemma 2.2 that
$res(h_1,h_0)\in U(R)$. Thus, $det(A)\in U(R)$, and so we can find
$y_i,z_j\in R$ such that $(*)$ and $(**)$ hold. In other wards, we
have $H_0$ and $H_1$ such that $H(g,t)=H_0H_1$. Obviously,
$H_0(g,0)=y_0+(\alpha_0-y_0)g$. As
$y_0+(\alpha_0-y_0)=\alpha_0=h_0(0)\in U(R)$, it follows by Lemma
3.5 that $H_0(g,0)\in U(RG)$, i.e., $H_0\in {\Bbb S}_0$. Further,
$H_1(g,1)=1+\sum\limits_{i=1}^{s-1}\big(z_i+(\beta_i-z_i)g\big)=1+\sum\limits_{i=1}^{s-1}z_i+
\big(\sum\limits_{i=1}^{s-1}(\beta_i-z_i)\big)g.$ It is easy to
check that $1+\sum\limits_{i=1}^{s-1}z_i+
\big(\sum\limits_{i=1}^{s-1}(\beta_i-z_i)\big)=1+\sum\limits_{i=1}^{s-1}\beta_i=h_1(1)\in
U(R)$. In view of Lemma 3.5, $H_1\in {\Bbb S}_1$. Clearly,
$\varphi(g):=res(H_0,H_1)\in RG$. As $\varphi(1)=res\big(H_0(1,t),
H_1(1,t)\big)=res(h_0,h_1)\in U(R)$. By using Lemma 3.5 again,
$\varphi(g)\in U(RG)$, i.e., $res(H_0,H_1)\in U(RG)$. In light of
Lemma 3.5, we get $(H_0,H_1)=1$. Therefore, $A(g)\in
M_n\big(RG\big)$ is strongly clean, as required.

$(2)\Rightarrow (1)$ Let $\psi: RG\to R, a+bg\mapsto a+b$. Then we
get a corresponding ring morphism $\mu: M_n(RG)\to M_n(R),
\big(a_{ij}(g)\big)\mapsto \big(\psi(a_{ij}(g))\big)$. As $A(g)$
is strongly clean, we can find an idempotent $E\in M_n(RG)$ such
that $A(g)-E\in GL_n(RG)$ and $EA=AE$. Applying $\mu$, we get
$A(1)-\mu(E)\in GL_n(R)$, where $\mu(E)\in M_n(R)$ is an
idempotent, hence the result. \hfill$\Box$

\vskip4mm \hspace{-1.8em} {\bf Example 3.7.}\ \ Let $S=\{
0,1,a,b\}$ be a set. Define operations by the following tables:
$$\begin{tabular}{cccc}
\begin{tabular}{c|c}
+ & \begin{tabular}{cccc} 0&1&a&b
\end{tabular}\\
\hline
\begin{tabular}{c}
0\\
1\\
a\\
b
\end{tabular}&
\begin{tabular}{cccc}
0&1&a&b\\
1&0&b&a\\
a&b&0&1\\
b&a&1&0
\end{tabular}
\end{tabular}
&,&&
\begin{tabular}{c|c}
$\times$&
\begin{tabular}{cccc} 0&1&a&b
\end{tabular}\\
\hline
\begin{tabular}{c}
0\\
1\\
a\\
b
\end{tabular}&
\begin{tabular}{cccc}
0&0&0&0\\
0&1&a&b\\
0&a&b&1\\
0&b&1&a
\end{tabular}
\end{tabular}
\end{tabular}.$$
Then $S$ is a finite field with $|~S~|=4$. Let $$R=\{
s_1+s_2z~|~s_1,s_2\in S, z~\mbox{is an indeterminant
satisfying}~z^2=0\}.$$ Then $R$ is a commutative local ring with
$charR=2$. We claim that $$A(x)=\left(
\begin{array}{cc}
\overline{az}&\overline{z+x}\\
\overline{1+x}&\overline{b+zx}
\end{array}
\right)\in M_2\big(R[x]/(x^2-1)\big)$$ is strongly clean. Clearly,
$A(1)=\left(
\begin{array}{cc}
az&1+z\\
0&1+bz
\end{array}
\right)\in M_2(R)$. As $\chi\big(A(1)\big)$ has a root $az\in
J(R)$ and a root $1+bz\in 1+J(R)$, $A(1)$ is strongly clean, and
we are through by Theorem 3.6.

\vskip30mm \bc {\Large\bf Acknowledgements}\ec \vskip4mm \no This
research was supported by the Natural Science Foundation of
Zhejiang Province (LY13A0 10019) and the Scientific and
Technological Research Council of Turkey (2221 Visiting Scientists
Fellowship Programme).

\end{document}